\documentclass{kms-b}
 \usepackage{amsmath}
 \usepackage{amsfonts}
 \usepackage{amssymb}
\usepackage{latexsym}

\usepackage{graphicx}
\allowdisplaybreaks

\theoremstyle{plain}
\newtheorem{theorem}{Theorem}[section]

\newtheorem{lemma}[theorem]{Lemma}

\theoremstyle{definition}

\theoremstyle{remark}
\newtheorem{remark}[theorem]{Remark}

\begin{document}

\title[Sobolev  trace inequality]
{Sobolev  trace inequality on $W^{s, q}( {\mathbb R}^{n} )$}

\author[Hee Chul Pak]{Hee Chul Pak}
\address{Hee Chul Pak \\
         Department of Mathematics \\ Dankook University  \\ Cheonan 330-714, South Korea}
\email{hpak@dankook.ac.kr}


\begin{abstract}
Sobolev trace inequalities on nonhomogeneous fractional Sobolev spaces
are established.
\end{abstract}

\maketitle

%%%%%%%%%%%%%%%%%%%%%%%%%%%%%%%%%%%%%%%%%%%%%%%%%%%%%%%%%%%%%%%%
\section{Trace inequalities on fractional Sobolev spaces}
%%%%%%%%%%%%%%%%%%%%%%%%%%%%%%%%%%%%%%%%%%%%%%%%%%%%%%%%%%%%%%%%

Sobolev trace inequality on nonhomogeneous
Sobolev spaces $H^s({\mathbb R}^{n})$ is given by:
for $1 \leq m<n$ and $s>\frac{m}{2}$,
\begin{align}
\| \tau u \|_{H^{s-\frac{m}{2}}({\mathbb R}^{n-m})}
 \leq  C_{s, m} \|u \|_{H^s({\mathbb R}^{n})},
         \label{NHS}
\end{align}
where $\tau u \in H^{s-\frac{m}{2}}({\mathbb R}^{n-m})$ is
the trace of $u \in H^s({\mathbb R}^{n})$
restricted to  the $(n-m)$-dimensional subspace
$
\{ x \in {\mathbb R}^{n} :
 x = (x_1, x_2, \cdots, x_{n-m}, 0, 0, \cdots, 0) \}$
 of $\mathbb{R}^n$.
The best constant for this inequality is presented by H. Pak snd Y. Park\cite{Pak} as
\[
C_{s, m}=
\left(
\frac{\Gamma(s-\frac{m}{2})}{(4\pi)^{\frac{m}{2}}\Gamma(s)}
\right)^{1/2}
\]
together with various forms of extremal functions.
A homogeneous version of trace inequality (\ref{NHS}) is
obtained by A. Einav and M. Loss  with the same sharp constant\cite{Loss}.
Even though the trace of a given function is addressed its fundamental importance
in the theory of boundary value problems of partial differential equations,
the continuity of trace operators
has not been reported yet even on general fractional Sobolev spaces
$W^{s, q}( {\mathbb R}^{n} )$.

This paper establishes fractional Sobolev  trace inequalities on the
nonhomogeneous fractional
Sobolev spaces $W^{s, q}( {\mathbb R}^{n} )$, $1< q \leq 2$.
The main result can be summarized as follows:

\begin{theorem}[Fractional Sobolev  trace inequality on $W^{s, q}( {\mathbb R}^{n} )$]
 \label{main}
Let $p, q$ be extended real numbers of
$\frac{1}{p} + \frac{1}{q} = 1$, $1 \leq  q \leq 2$
and let $s, t$ be real numbers with
\begin{align}
s -n \left(
               \frac{1}{q} - \frac{1}{p}
         \right) > t \geq \frac{m}{p}.
 \label{ind-cond}
\end{align}
Then
for  $u \in W^{s, q}( {\mathbb R}^{n} )$ with
the trace $\tau u$ on $\mathbb{R}^{n-m}$,
we have
\begin{align}
\| \tau u \|_{W^{t-\frac{m}{p}, p}({\mathbb R}^{n-m})} \!
 \leq  C_{q, s, t, m}
\|u \|_{W^{s, q}({\mathbb R}^{n})},  \label{STE}
\end{align}
where
\begin{align}
C_{q, s, t, m} =
\frac{\pi^{ \frac{n-m }{2q} -\frac{n }{2p}  }    }
{2^{\frac{m}{2}}}
    {q^{ { 2n-m \over 2q }  }  \over  p^{{2n -m \over 2p} } }
  \left[
     \frac{\Gamma(\frac{sq-m}{2} )}{\Gamma(\frac{sq}{2})}
  \right]^{{1 \over q}}
  \left[
     \frac{\Gamma(\frac{(s-t)p}{2(p-2)} - \frac{n}{2})}
          {\Gamma( \frac{(s-t)p}{2(p-2)} - \frac{m}{2} )}
  \right]^{\frac{1}{q} - \frac{1}{p} }.
\label{constant}
\end{align}
\end{theorem}  %%===============================================

The main difficulty of the proof of the theorem arises
from the absence
of the isometry, and even from the lack of the continuity
of the Fourier transform  on $L^p(\mathbb{R}^n)$ with
$p \neq 2$.

The classical  Sobolev trace inequalities on $\mathbb{R}^{n+1}$ are given
by
\begin{align}
\| \tau u \|_{L^p({\mathbb R}^{n})}
 \leq  A_{q, n}
\|\nabla u \|_{L^q({\mathbb R}^{n+1}) }  \label{classical-NHS}
\end{align}
with
\[
1 -n \left(
               \frac{1}{q} - \frac{1}{p}
         \right) =  \frac{1}{q},
\]
for some  positive constant $A_{q, n}$ which is independent of the function $u$.
It can be noticed that the range $L^p(\mathbb{R}^n)$ of the trace map is {\it too big}
to constitute a proper container of all the traces.
Theorem \ref{main} illustrates that
all the traces of the functions in $W^{1,  q}({\mathbb R}^{n+1})$
are  included at least in the  space
\[
\bigcap W^{t-\frac{1}{p},  p}({\mathbb R}^{n}),
\]
where the intersection is taken over all indices $t$ satisfying  (\ref{ind-cond}).

The best constant of (\ref{classical-NHS}) is still open  except for the cases $q=1$ and $q=2$.
A conjectured  extremal function is the function of the form
\[
\tau u(x) = {1 \over (1 +|x|^2)^{n+1 - p \over 2(p-1)} }.
\]
J. Escobar  first identified the best constant for the case $q=2$ of (\ref{classical-NHS})
by exploiting the conformal
 invariance of this inequality and using characteristics of an Einstein metric\cite{Escobar}.
W. Beckner  independently achieved the sharp constant  by inverting the inequality to
 a fractional integral on the dual space and using
 the sharp Hardy-Littlewood-Sobolev inequality\cite{Beckner2}.
The limiting case $p=1$ is investigated by Y. Park\cite{Park}.
We present an upper bound for the sharp constant of (\ref{STE}) and observe
that $C_{q, s, t, m}$ blows up to infinity
as $t$ approaches to $s -n \left(
               \frac{1}{q} - \frac{1}{p}
         \right)$.
We also note that
$C_{q, s, t, m}$ converges to $0$ as $s$ goes to infinity.

Some basic notations are listed.
The fractional Sobolev spaces
$W^{s,  q}({\mathbb R}^{n})$ of functions with $s \in \mathbb R$
are defined as
\[
W^{s,  q}({\mathbb R}^{n})
:=
        \left\{
              u \in {\mathcal S}'({\mathbb R}^{n}) :
              {\mathcal F}_n^{-1}
              \left(
              (1+ |\xi|^2)^{s/2} \widehat{u}
              \right) \in L^q({\mathbb R}^{n})
        \right\},
\]
where ${\mathcal S}'({\mathbb R}^{n})$ is the set of all tempered
distributions on ${\mathbb R}^{n}$ and the Fourier transform
$\widehat{u}= {\mathcal F}_n(u)$ on ${\mathbb R}^{n}$ of the function $u
\in {\mathcal S}({\mathbb R}^{n})$ is defined by
\[
\widehat{u}(\xi) = {\mathcal F}_n(u)(\xi) =  \frac{1}{(2\pi)^{n/2}}
 \int_{{\mathbb R}^{n}} u(x) e^{-i x \cdot \xi}  \, dx.
\]
The nonhomogeneous Sobolev space $W^{s, q}({\mathbb R}^{n})$  is equipped with
the norm
\[
\|u \|_{W^{s, q}} := \left(
 \int_{{\mathbb R}^{n}}
       \left| {\mathcal F}_n^{-1}
              \left(
              (1+ |\xi|^2)^{s/2} \widehat{u}
              \right)(x)
       \right|^q \, dx
 \right)^{1/q}.
\]

%%%%%%%%%%%%%%%%%%%%%%%%%%%%%%%%%%%%%%%%%%%%%%%%%%%%%%%%%%%%%%%%
\section{Proof of Theorem \ref{main} }
%%%%%%%%%%%%%%%%%%%%%%%%%%%%%%%%%%%%%%%%%%%%%%%%%%%%%%%%%%%%%%%%

The nonhomogeneous Sobolev space $W^{s, q}({\mathbb R}^{n})$
is the completion of
Schwartz class ${\mathcal S}({\mathbb R}^{n})$.
Hence, by the continuous extension argument,  it suffices to show that
for $u \in S({\mathbb R}^{n})$, we have
\begin{align}
\| \tau u \|_{W^{t-\frac{m}{p}, p}({\mathbb R}^{n-m})}
 \leq  \!
C_{q, s, t, m}
\|u \|_{W^{s, q}({\mathbb R}^{n})},
         \label{s-STE}
\end{align}
where $C_{q, s, t, m}$ is the constant defined at (\ref{constant}).
Here  the trace $\tau u$ of $u$ is defined by
$\tau u(x') = u(x', 0, \cdots, 0)$ for $x' \in {\mathbb R}^{n-m}$.

To accomplish it,
take $u \in {\mathcal S}({\mathbb R}^{n})$ and set $f := \tau u$ to have
\begin{align*}
\widehat{f}(\xi') = {\mathcal F}_{n-m}(f)(\xi ')
&=
\frac{1}{(2\pi)^{ \frac{n-m}{2}  }}
\int_{{\mathbb R}^{n-m}} f(x') e^{-i x'\cdot \xi'} dx' \\
&=
\frac{1}{(2\pi)^{\frac{n-m}{2} }}
\int_{{\mathbb R}^{n-m}} u(x', 0) e^{-i x'\cdot \xi'} dx'
\end{align*}
for $\xi' \in {\mathbb R}^{n-m}$.
Apply the Fourier inversion formula in the $\xi''$-variable to get
\[
u(x', 0) = \frac{1}{(2\pi)^{m/2}}
            \int_{{\mathbb R^m}}
                         {\mathcal F}_m(u(x', \cdot \,)) (\xi'')
                         %e^{i 0 \cdot \xi ''}
                         d \xi''
\]
for $x' \in {\mathbb R}^{n-m}$,
where ${\mathcal F}_m$ represents the Fourier transform with respect to
$x''$-variable for $(x', x'' ) \in {\mathbb R}^{n-m} \times \mathbb{R}^m$.
Then Fubini's theorem gives
\begin{align}
\widehat{f}(\xi')
&=   \frac{1}{(2\pi)^{n/2}}
          \int_{{\mathbb R}^{n-m}}
          \left(
              \int_{{\mathbb R^m}}
               {\mathcal F}_m (u(x',  \cdot \,)) (\xi '') \; d \xi''
          \right)
               e^{-i x'\cdot \;  \xi'} dx'
     \nonumber \\
&=   \frac{1}{(2\pi)^{n/2}}
    \int_{{\mathbb R^m}}  \!\!
    \int_{{\mathbb R}^{n-m}}  \!\!\!
     \left(
    \frac{1}{(2\pi)^{m/2}}
    \int_{{\mathbb R^m}}
         u(x',x'') e^{-i x'' \cdot \xi''} dx'' \right)  \!
    e^{-i x'\cdot \; \xi'} dx'
    d \xi'' \nonumber \\
&=  \frac{1}{(2\pi)^{m/2}}
     \int_{{\mathbb R^m}}
            \widehat{u}(\xi', \xi'') \,
            d \xi'',   \label{tr-A}
\end{align}
where $(\xi', \xi ''), (x', x '') \in {\mathbb R}^{n-m}\times{\mathbb R}^{m}$.
Let $\phi$ be a function  in $S({\mathbb R}^{n-m})$ with
$\| \phi  \|_{L^q({\mathbb R}^{n-m})} = 1$.
Multiply both sides of
(\ref{tr-A}) by $(1+ |\xi'|^2)^{{t \over 2} - \frac{m}{2p}  }  \widehat{\overline{\phi}}(- \xi')   $
and integrate with respect to $\xi'$ to get
\begin{align}
\int_{{\mathbb R^{n-m}}}
&\widehat{f}(\xi')(1+ |\xi'|^2)^{{t \over 2} - \frac{m}{2p}  }
    \widehat{\overline{\phi}}(- \xi')  d \xi'    \nonumber \\
& = \! \frac{1}{(2\pi)^{{m \over 2} }}
   \int_{{\mathbb R^{n}}} \!\!
        \widehat{u}(\xi', \xi'') (1+ |\xi'|^2)^{{t \over 2} - \frac{m}{2p}  }
          \widehat{\overline{\phi}}(- \xi') d \xi
  \label{tr-B}
\end{align}
for $\xi = (\xi', \xi'') \in {\mathbb R}^{n-m} \times {\mathbb R}^{m}$.
Here the bar $\overline{z}$ indicates the complex conjugate of $z$.
The left hand side of (\ref{tr-B}) becomes
\begin{align}
\int_{{\mathbb R^{n-m}}}  \!\!\!\!
  \widehat{f}(\xi')
&(1+ |\xi'|^2)^{{t \over 2} - \frac{m}{2p}  }
\widehat{\overline{\phi}}(- \xi')   d \xi'
\nonumber \\
& =
 \int_{{\mathbb R^{n-m}}}  \!\!
\left[
\widehat{f}(\xi')(1+ |\xi'|^2)^{{t \over 2} - \frac{m}{2p}  }
\right]^{\! \vee} \!\!\! (x') \,
  \overline{\phi (x')}  \, d x',
   \label{L-side}
\end{align}
where  $g^{\! \vee} $ represents  the Fourier inversion ${\mathcal F}_n^{-1}(g)$
of the function $g$.
H\"{o}lder's inequality on the right hand side of (\ref{tr-B}) yields that
\begin{align}
&\left|
\frac{1}{(2\pi)^{{m \over 2} }}  \!\!\!
   \int_{{\mathbb R^{n}}} \!\!
        \widehat{u}(\xi) (1+ |\xi|^2)^{\frac{s}{2} }  {\overline{\phi}}^{ \vee}(\xi')
        \frac{(1+ |\xi'|^2)^{{t \over 2} - \frac{m}{2p}  } }
             {(1+ |\xi|^2)^{\frac{s}{2} } }
          d \xi
\right| \nonumber \\
&\leq
\frac{1}{(2\pi)^{{m \over 2} }}  \!\!\!
\left(
\int_{{\mathbb R^n}}  \!\!
            \left|
               \widehat{u}(\xi) (1+ |\xi|^2)^{\frac{s}{2}}
            \right|^p d \xi
\right)^{  \!\! {1 \over p}}  \!\!\!
\left(
\int_{{\mathbb R^n}}  \!\!
            | {\overline{\phi}}^{ \vee}(\xi') |^{q}
        \frac{(1+ |\xi'|^2)^{{tq \over 2} - \frac{mq}{2p}  } }
             {(1+ |\xi|^2)^{\frac{sq}{2} } }
         d \xi
\right)^{ \!\! {1 \over q}}   \label{HolderIneq_p}
\end{align}
with $\frac{1}{p} + \frac{1}{q} = 1$.

\begin{lemma} \label{lemm-1}
For $\alpha > \frac{m}{2}$, we have
\begin{align}
 \int_{{\mathbb R}^{m}}\frac{(1+ |\xi'|^2)^{\alpha- \frac{m}{2} }}
 {(1+ |\xi'|^2 + |\xi''|^2)^{\alpha}} \, d \xi''
&= \pi^{\frac{m}{2}}\:
\frac{\Gamma(\alpha -\frac{m}{2})}{\Gamma(\alpha)}.
\label{Holder_Ineq_p-A}
\end{align}
\end{lemma}

\proof
A direct computation reveals that
for $\xi = (\xi', \xi_n ) \in {\mathbb R}^{n-1} \times {\mathbb R}$,
\begin{align*}
 \int_{{\mathbb R}} \frac{1}{(1+ |\xi'|^2+\xi_n^2)^{\alpha}} d \xi_n
&=  \int_{-\frac{\pi}{2}}^{\frac{\pi}{2}}
     \frac{1}{(1+ |\xi'|^2 )^{\alpha - \frac{1}{2}}}
      (1+ \tan^2 \theta)^{-\alpha  +1} d \theta \\
&=   \frac{1}{(1+ |\xi'|^2 )^{\alpha - \frac{1}{2}}}
\left( 2  \int_0^{\frac{\pi}{2}} \cos^{2\alpha -2}\theta  d \theta \right) \\
&=  \sqrt{\pi} \;\frac{\Gamma(\alpha - \frac{1}{2})}
         {\Gamma(\alpha)}
    \frac{1}{(1+ |\xi'|^2 )^{\alpha - \frac{1}{2}}}.
\end{align*}
An induction argument gives that for fixed $\xi' \in {\mathbb R}^{n-m}$,
\begin{align*}
 \int_{{\mathbb R}^{m}}\frac{1}{(1+ |\xi'|^2 + |\xi''|^2)^{\alpha}} \, d \xi''
&= \pi^{\frac{m}{2}}\:
\frac{\Gamma(\alpha  -\frac{m}{2})}{\Gamma(\alpha)}
\frac{1}{(1+ |\xi'|^2)^{\alpha - \frac{m}{2} }}.
\end{align*}
This implies the identity (\ref{Holder_Ineq_p-A}).
\hfill$\Box$\par
  %%****************************************************************

\bigskip

Lemma \ref{lemm-1}
and H\"{o}lder's inequality   imply that
\begin{align}
\int_{{\mathbb R^n}}
  &| {\overline{\phi}}^{ \vee}(\xi') |^{q}
            \frac{(1+ |\xi'|^2)^{{tq \over 2} - \frac{mq}{2p}  } }
             {(1+ |\xi|^2)^{\frac{sq}{2} } }
         d \xi   \nonumber  \\
&=
\pi^{\frac{m}{2}}\:
\frac{
           \Gamma(\frac{sq}{2}-\frac{m}{2})}{\Gamma(\frac{sq}{2})}
\int_{{\mathbb R^{n-m}}}
   \frac{| {\overline{\phi}}^{ \vee}(\xi') |^{q} }
             {(1+ |\xi'|^2)^{\frac{(s-t)q}{2} + \frac{m(q-2)}{2}} }
         d \xi'   \nonumber  \\
& \leq
\pi^{\frac{m}{2}}\:
\frac{
           \Gamma(\frac{sq}{2}-\frac{m}{2})}{\Gamma(\frac{sq}{2})}
\left(
\int_{{\mathbb R^{n-m}}}
   | {\overline{\phi}}^{ \vee}(\xi') |^{p}
         d \xi'
\right)^{q \over p}
 \nonumber  \\
& \qquad \qquad \times
\left(
\int_{{\mathbb R^{n-m}}}
   \frac{1 }
             {(1+ |\xi'|^2)^{\frac{(s-t)p}{2(p-2)} - \frac{m}{2      } } }
         d \xi'
\right)^{{q \over p} (p-2)}    \nonumber  \\
&=
\pi^{\frac{m}{2}}\:
\frac{\Gamma(\frac{sq}{2}-\frac{m}{2})}{\Gamma(\frac{sq}{2})}
\left[
\frac{\pi^{\frac{n-m}{2}} \Gamma(\frac{(s-t)p}{2(p-2)}  - \frac{n}{2})}
     {\Gamma(\frac{(s-t)p}{2(p-2)}  - \frac{m}{2} )}
\right]^{{q \over p} (p-2)}
\!\!\!\!
\left(
\int_{{\mathbb R^{n-m}}}
   | {\overline{\phi}}^{ \vee}(\xi') |^{p}
         d \xi'
\right)^{q \over p}.
            \label{Comp-1}
\end{align}
By virtue of the Babenko-Beckner's inequality \cite{Babenko, Beckner1},
we have
\begin{align}
\left(
\int_{{\mathbb R^{n-m}}}
   | {\overline{\phi}}^{ \vee}(\xi') |^{p}
         d \xi'
\right)^{1 \over p}
\leq
{q^{ { n-m \over 2q }  }  \over  p^{{n -m \over 2p} } }
\left(
\int_{{\mathbb R^{n-m}}}
   | \overline{\phi(\xi') } |^{q}
         d \xi'
\right)^{1 \over q}
=
{q^{ { n-m \over 2q }  }  \over  p^{{n -m \over 2p} } }.
            \label{test-fn}
\end{align}
and
\begin{align}
\left( \int_{{\mathbb R^n}}
            \left|  \widehat{u}(\xi) (1+ |\xi|^2)^{\frac{s}{2}} \right|^p \, d \xi
            \right)^{1/p}   \!\!\! % \nonumber \\
&=
\left( \int_{{\mathbb R^n}}
            \left|
              {\mathcal F}_{n} \circ {\mathcal F}_{n}^{-1}
              \left(   \widehat{u}(\xi) (1+ |\xi|^2)^{\frac{s}{2}}  \right)
            \right|^p \, d \xi
            \right)^{1/p}  \nonumber \\
&\leq
 {q^{{n \over 2q} } \over  p^{ { n \over 2p }  } }
\left( \int_{{\mathbb R^n}}
            \left|
               \left(  \widehat{u}(\xi) (1+ |\xi|^2)^{\frac{s}{2}}   \right)^{\! \vee}(x)
            \right|^q  \, d x
            \right)^{1/q}.
               \label{Holder_Ineq_p-2}
\end{align}
Collecting the estimates (\ref{tr-B}), (\ref{L-side}),  (\ref{HolderIneq_p}),
 (\ref{Comp-1}), (\ref{test-fn})
and (\ref{Holder_Ineq_p-2}),
we establish the inequality (\ref{s-STE}).
The proof is now completed.
\hfill$\Box$\par
%%****************************************************************

\bigskip

\begin{remark}
Let $p, q$ be extended real numbers of
$\frac{1}{p} + \frac{1}{q} = 1$, $1 \leq  q \leq 2$
and let $s, t$ be real numbers with
$
t - n \left(
               \frac{1}{q} - \frac{1}{p}
         \right)
> s > \frac{t}{2}+ \frac{m}{p}
$.
Then
for any $g \in W^{t -\frac{m}{q}, q}({\mathbb R}^{n-m})$,
there is a function $u \in W^{s, p}({\mathbb R}^{n}) $ such that
\[
\tau u = g.
\]
In fact, for $g \in W^{t -\frac{m}{q}, q}({\mathbb R}^{n-m})$,
we consider the function
$$
u(x)
:=
 2^{m/2}
\frac{\Gamma(s )}{\Gamma(s-\frac{m}{2})}
 \left( \widehat{g}(\xi')
\frac{
(1+ |\xi'|^2)^{s -\frac{m}{2} }  }
     {    (1+ | \xi  |^2)^{ s }      }
 \right)^{\! \vee}(x).
$$
Then by the identity (\ref{tr-A}) together with Lemma \ref{lemm-1}, we observe that
\begin{align*}
\widehat{\tau u}(\xi') =  \frac{1}{(2\pi)^{m/2}}
     \int_{{\mathbb R^m}}
            \widehat{u}(\xi', \xi'') \,
            d \xi''
     =    \widehat{g}(\xi').
\end{align*}
In order to demonstrate that $u$ belongs to $W^{s, p}({\mathbb R}^{n}) $,
we first note that
\begin{lemma} \label{lemm-2}
For $\alpha > \frac{m}{p}$ and $\beta \in \mathbb{R}$, we have
\begin{align*}
\int_{{\mathbb R}^{n}}
\left|
\widehat{g}(\xi')
 \frac{(1+ |\xi'|^2)^{\beta - \frac{m}{2} }}
 {(1+ |\xi|^2)^{\alpha}}
\right|^p
\, d \xi
&=
C
\int_{{\mathbb R}^{n-m}}
\left|
\, \widehat{g}(\xi')
 (1+ |\xi'|^2)^{\beta-\alpha - \frac{m}{2q} }
 \right|^p
\, d \xi'
\end{align*}
for some constant $C$.
\end{lemma}
\proof
The result follows from a direct computation and Lemma \ref{lemm-1}.
\hfill$\Box$\par
  %%****************************************************************

\bigskip

The same arguments used in the proof of Theorem \ref{main} yields: for any
$\| \phi  \|_{L^q(\mathbb{R}^n)} = 1$
\begin{align*}
\left|
\int_{{\mathbb R^n}} \!\!\!\!
           \left(
              \widehat{u}(\xi)  (1+ |\xi|^2)^{ \frac{s}{2} }
           \right)^{ \vee} \!
            \overline{\phi(x) }
\, d x
\right|  \!
&=  C_1
\left|
\int_{{\mathbb R^n}}
              \widehat{g}(\xi')
 \frac{(1+ |\xi'|^2)^{s - \frac{m}{2} }}
 {(1+ |\xi|^2)^{ \frac{s}{2} } } \,
            \overline{\phi}^{ \vee}(\xi)
\, d \xi
\right|   \\
&\lesssim
\left(
\int_{{\mathbb R^n}}
          \left|
           \widehat{g}(\xi')
 \frac{(1+ |\xi'|^2)^{s - \frac{m}{2} }}
 {(1+ |\xi|^2)^{ s- \frac{t  }{2}   } } \,
      \right|^p
\, d \xi
\right)^{1/p}      \\
&\qquad \qquad \times
\left(
\int_{{\mathbb R^n}}
\left|
 \frac{ \overline{\phi}^{ \vee}(\xi) }
 {(1+ |\xi|^2)^{ \frac{t - s}{2} } }
\right|^q
\, d \xi
\right)^{1/q}   \\
&= C_2
\left(
\int_{{\mathbb R^{n-m}}}
            \left|
             \left(
                  \widehat{g}(\xi') (1+ |\xi'|^2 )^{ \frac{t}{2} - \frac{m}{2q}  }
             \right)
            \right|^p  \, d \xi'
\right)^{1/p}   \\
&\lesssim  \!
\left(
\int_{{\mathbb R^{n-m}}}  \!\!
            \left|
             \left(
               \widehat{g}(\xi') (1+ |\xi'|^2 )^{ {t  \over 2} - \frac{m}{2q}   }
             \right)^{\! \vee} \!\!\! (x')
            \right|^q  d x'
\right)^{1/q},
\end{align*}
for some positive constants $C_1$ and $C_2$.
Hence  we can see that  $u$ belongs to $W^{s, p}({\mathbb R}^{n}) $ and $\tau u = g$.
\hfill$\Box$\par
%%****************************************************************
\end{remark}

\vspace{\baselineskip}

\end{document}